\documentclass[12pt]{article}
\usepackage{amsmath,amsthm,amssymb,amscd}
\usepackage{latexsym}
\newtheorem{thm}{Theorem}[section]
 
 \newtheorem{lem}[thm]{Lemma}
 
 \theoremstyle{definition}
 
 \theoremstyle{remark}

 \numberwithin{equation}{section}

\title
{Normal curvature bounds along the mean curvature flow
\thanks{Partially supported by NSFC no.10671018.}}

\author{ Hong Huang$^{\star}$
  \\
\small School of Mathematical Sciences,
Beijing Normal University \\
\small Beijing 100875,
P.R. China\\
\small {\em E-mail address:} {\bf hhuang@bnu.edu.cn}\\
}
\date{}
\begin{document}
\maketitle
\begin{abstract}

Let $(M^n,g_0)$ and $(\bar{M}^{n+1},\bar{g})$ be
complete Riemannian manifolds with $|\bar{\nabla}^k\bar{Rm}|\le
\bar{C}$ for $k \le 2$, and suppose there is an isometric immersion
$F_0: M^n \rightarrow \bar{M}^{n+1}$ with bounded second fundamental form. Let $F_t: M^n \rightarrow
\bar{M}^{n+1}$ ($t\in [0,T]$) be a family of immersions evolving by
mean curvature flow with initial data $F_0$ and with uniformly
bounded second fundamental forms.
 We show that the supremum  and infimum of the
normal curvature of the immersions $F_t$ vary at a bounded rate. This
is an analogue of  a result of Rong and Kapovitch on Ricci flow.

{\bf Key words}: mean curvature flow, barrier function, maximum
principle

{\bf AMS2010 Classification}: 53C44
\end{abstract}
\maketitle

\section {Introduction}

 It is proved in Rong \cite{R}  that the supremum  and
infimum of the sectional curvature of a compact  manifold  vary at a
bounded rate under the Ricci flow. In a recent paper \cite{K} Kapovitch tried to extend this  result  to the noncompact case. 
In this short note, we will prove a
mean curvature flow analogue of their results. More precisely, we have
the following

 \begin{thm} \label{thm} Let $(M^n,g_0)$ and $(\bar{M}^{n+1},\bar{g})$ be
complete Riemannian manifolds with $|\bar{\nabla}^k\bar{Rm}|\le
\bar{C}$ for $k \le 2$, and suppose there is an isometric immersion
$F_0: M^n \rightarrow \bar{M}^{n+1}$ with bounded second fundamental form. Let $F_t: M^n \rightarrow
\bar{M}^{n+1}$ ($t\in [0,T]$) be a family of immersions evolving by
mean curvature flow with initial data $F_0$ and with uniformly
bounded second fundamental forms $|A|\le \bar{C}$. Then there exists
a constant $C$ (depending only on $\bar{C}, T$ and $n$ ) such that $ inf
\kappa_0-Ct \leq \kappa_t \leq sup \kappa_0 +Ct$, where $ \kappa_t$
is the normal curvature function of the immersion $F_t$.

\end{thm}

Recall that the short time existence of the mean curvature flow  in our situation is well-known, cf., for example, Theorem 4.2 in Ecker-Huisken \cite{EH} and
Proposition 10.3 in Zhu \cite{Z}. If $M$ is compact, the uniqueness of the codimension one mean curvature flow is also well-known.  If $M$ is not compact the uniqueness is established by Chen and Yin \cite{CY} under  an extra condition on the injectivity radius of $\overline{M}$ (Chen and Yin's result also applies to the higher codimension case).

Moreover from   Chen and Yin \cite{CY}(in particular, Corollary 3.3 there) we know that under the condition of our theorem the following
estimates hold on $M^n \times [0,T]$:

$|Rm_{g_t}|\leq C$, $|\nabla Rm|_{g_t}\leq \frac{C}{t^{1/2}}$  and $|\frac{\partial g_t}{\partial t}|\leq C$,

\noindent where $g_t$ is the induced metric on $M$ from the immersion
$F_t$, and  $C$ is a constant
depending only on  $\bar{C}, T$ and $n$. 

  In the next section we will use the barrier function technique as in Section 2.3 of Cao and Zhu \cite{CZ} (instead of the cut-off function technique in  Kapovitch \cite{K}) combined with a maximum principle argument to prove our result.

\section {Proof of Theorem}

  For simplicity we only consider
 the case $\overline{M}^{n+1}=R^{n+1}$, the general case can  be treated
 similarly, since in the general case one need only
add some lower order terms in the evolution equation of the second
fundamental form, which hardly affect the  proof below.

First we prove the existence of a nice barrier function similarly as in the Ricci flow literatures (cf., for example, Lemma 2.1.1 in Cao and Zhu \cite{CZ}, and Lemma 12.5 in Chow et al \cite{Ch}).

\begin{lem} \label{lem} Let $(M^n,g_0)$ and $(\bar{M}^{n+p},\bar{g})$ be
complete Riemannian manifolds with $|\bar{\nabla}^k\bar{Rm}|\le
\bar{C}$ for $k \le 2$, and suppose there is an isometric immersion
$F_0: M^n \rightarrow \bar{M}^{n+p}$ with bounded second fundamental form. Let $F_t: M^n \rightarrow
\bar{M}^{n+p}$ ($t\in [0,T]$) be a family of immersions evolving by
mean curvature flow with initial data $F_0$ and with uniformly
bounded second fundamental forms $|A|\le \bar{C}$. Then there exists
a smooth function $f$ on $M$ satisfying $f \ge 1$, $f(x)\rightarrow
\infty$ as $d_{g_0}(x,x_0)\rightarrow \infty$ (for some fixed point
$x_0 \in M$), $|\nabla f|_{g_t} \le C$, and $|{\nabla}^2 f|_{g_t}
\le C$ for some constant $C$ (depending only on $\bar{C}, T,  n$ and
$p$), where $g_t=F_t^*\bar{g}$.
\end{lem}

\vspace *{0.4cm}

 Proof. The proof is along the lines of that of
Lemma 2.1.1 in Cao and Zhu \cite{CZ}, and Lemma 12.5 in Chow et al \cite{Ch},
with Shi's derivative estimate replaced by Corollary 3.3 in Chen and
Yin \cite{CY}.

\vspace *{0.4cm}

Now let $\kappa_t(x,v)=\frac{h_{ij}(x,t)v^iv^j}{|v|_{g_t}^2}$
be the normal curvature function (for the immersion $F_t$) in a direction $v\in T_xM$ at a point $x\in M$, where $A=(h_{ij})$ is the second fundamental form, $v^i$ are the components of $v$ in a coordinate system. Let
$\phi(t,x,v)=\kappa_t(x,v)-\varepsilon e^{Bt}f(x)$  as in the proof of Theorem 2.3.1 in Cao and Zhu \cite{CZ}, where $\varepsilon$ is a (small) positive constant, $B$ is a (large) constant to be chosen later, and $f$ is the barrier function given in Lemma \ref{lem}. Finally  we let
 $s(t)=sup \{\phi(t,x,v)\}$,
where $x$ runs over $ M$, and $v$ runs over nonzero vectors in
$T_xM$.

Given $t_0 \in [0,T]$, let $x_0 \in M$, and $v_0 \in T_{x_0}M$
 be such that $s(t_0)=\phi(t_0,x_0,v_0)$. We want to use Lemma 3.5 in Hamilton \cite{Ha} to estimate $\frac{ds}{dt}$, so we try to  show that
  $\frac{\partial \phi}{\partial t}(t_0,x_0,v_0)\leq C$.

Using the evolution equation for the second fundamental form

$$
  \frac{\partial}{\partial t}h_{ij}=\triangle h_{ij}-2Hh_{il}g^{lm}h_{mj}+|A|^2 h_{ij}
  =:\triangle h_{ij}+P_{ij}$$,

 \noindent(cf. Huisken  \cite{Hu}) we compute in a normal coordinate system at $x_0$ w.r.t
 $g_{t_0}$

$$\frac{\partial \phi}{\partial t}(t_0,x_0,v_0)
 =\triangle
h_{ij}(x_0,t_0)\frac{{v_0}^i{v_0}^j}{|v_0|_{g_{t_0}}^2}+
P_{ij}(x_0,t_0)\frac{{v_0}^i{v_0}^j}{|v_0|_{g_{t_0}}^2}
+ h_{ij}(x_0,t_0){v_0}^i{v_0}^j \frac{\partial}{\partial
t}|_{t_0}(\frac{1}{|v_0|_{g_t}^2})-\varepsilon Be^{Bt_0}f(x_0).$$

Now we extend the vector $v_0$ by parallel translation along
geodesics emanating radially out of $x_0$ w.r.t. $g_{t_0}$. Still
denote this vector field by $v_0$. Note 

$$\triangle (\frac
{h_{ij}(x,t_0){v_0}^i(x){v_0}^j(x)}{|v_0(x)|_{g_{t_0}}^2}-\varepsilon
e^{Bt_0}f(x)) \leq 0$$ 

\noindent at $x=x_0$, since the function
$\frac{h_{ij}(x,t_0){v_0}^i(x){v_0}^j(x)}{|v_0(x)|_{g_{t_0}}^2}-\varepsilon
e^{Bt_0}f(x)$ has a local maximum at $x_0$. As in Rong \cite{R} we have
$|\nabla ^2 {v_0}| (x_0, t_0)\leq C |v_0|_{g_{t_0}}$, so $\triangle
\frac{{v_0}^i{v_0}^j}{|v_0|_{g_{t_0}}^2}(x_0)$ is bounded (noting
$|v_0|_{g_{t_0}}$ does not depend on $x$). Clearly we  have
 $|\frac{\partial}{\partial
t}|_{t_0}|v_0|_{g_t}|\leq C |v_0|_{g_{t_0}}$. Then using Lemma \ref{lem} we have

$$\frac{\partial \phi}{\partial t}(t_0,x_0,v_0)\le C+\varepsilon
e^{Bt_0}(C-Bf(x_0)).$$

\noindent Choosing $B$ large enough such that $C-Bf(x_0)\le 0$, we get

$$\frac{\partial \phi}{\partial t}(t_0,x_0,v_0)\le C.$$

\noindent It follows from Lemma 3.5 in Hamilton \cite{Ha} that $\frac{ds}{dt} \le C$, and
$s(t)\le s(0)+Ct$, i.e.

$$sup \{\kappa_t(x,v)-\varepsilon e^{Bt}f(x)\} \le sup
\{\kappa_0(x,v)-\varepsilon f(x)\}+Ct,$$ 

\noindent which implies

$$ \kappa_t(x,v)-\varepsilon e^{Bt}f(x) \le sup
\{\kappa_0(.,.)\}+Ct.$$

\noindent Letting $\varepsilon \rightarrow 0$, we get

$$ \kappa_t(x,v) \le sup \{\kappa_0(.,.)\}+Ct,$$

\noindent  so 

$$sup
\{\kappa_t(x,v)\} \le sup \{\kappa_0(x,v)\}+Ct.$$

The argument for $inf
\kappa_t$ is similar.

\hspace *{0.4cm}

{\bf Remark}  Kapovitch \cite{K} tried to prove a similar result in Ricci
flow using a cut-off function argument. But it seems to me that his
argument may need justifications. That is it is not clear why
${\overline{A}_z}'(t)\leq C$ implies ${\overline{A}}'(t)\leq C$
there. Note that one cannot apply directly Lemma 3.5 in Hamilton \cite{Ha} in that situation, since $z$ runs over the whole noncompact manifold
$M$ rather than a compact subset of it. We use a barrier function argument instead, which can also be applied to the Ricci flow case and give an alternative proof of Kapovitch's result.


\hspace *{0.4cm}


\begin{thebibliography}{99}



\bibitem{CZ} H.-D. Cao, X.-P. Zhu, A complete proof of the
Poincar$\acute{e}$ and geometrization conjectures- application of
the Hamilton-Perelman theory of the Ricci flow, Asian J. Math. 10
(2006), 165-492.








\bibitem{CY}B.-L. Chen, L. Yin, Uniqueness and pseudolocality theorems of the mean curvature
flow, Comm. Anal. Geom. 15 (2007), no.3, 435-490.


\bibitem{Ch} B. Chow, S.-C. Chu, D. Glickenstein, C. Guenther, J. Isenberg, T. Ivey,
D. Knopf, P. Lu, F. Luo, and  L. Ni, The Ricci flow: techniques and
applications. Part II. Analytic aspects. Mathematical Surveys and
Monographs, 144. American Mathematical Society, Providence, RI,
2008.





\bibitem{EH} K. Ecker, G. Huisken, Interior estimates for hypersurfaces moving by mean curvature,
Invent. Math. 105 (1991), 547-569.

\bibitem{Ha} R. Hamilton, Four-manifolds with positive curvature operator, J. Diff. Geom. 24 (1986), 153-179.


\bibitem{Hu} G. Huisken, Flow by mean curvature of convex
surfaces into spheres, J. Diff. Geom. 20 (1984), no. 1, 237-266.



\bibitem{K} V. Kapovitch, Curvature bounds via Ricci smoothing, Illinois J. Math. 49 (2005),
no.1, 259-263.



\bibitem{R} X. Rong, On the fundamental groups of manifolds of positive sectional curvature,
Ann. Math. 143 (1996), no.2, 397-411.

\bibitem{Z} X.P. Zhu, Lectures on mean curvature flow, AMS/IP Studies in Advanced Mathematics, vol 32, Amer. Math. Soc. and Internat. Press. 2002.


\end{thebibliography}
\end{document}